\documentclass{ifacconf}
\usepackage{amsmath, amssymb}

\counterwithin*{section}{part}
\usepackage{cleveref}

\usepackage{graphicx,subfigure}      
\usepackage[numbers]{natbib}        
\usepackage{tikz}
\usetikzlibrary{arrows.meta}
\usetikzlibrary{cd,calc}
\tikzset{>=latex}
\graphicspath{{./Figs/}}

\usepackage{dsfont}

\definecolor{myblue}{RGB}{87, 163, 213}
\definecolor{myorange}{RGB}{226, 112, 46}
\definecolor{myred}{RGB}{218, 0, 0}
\definecolor{mygreen}{RGB}{129, 214, 83}
\definecolor{mygray}{gray}{0.9}

\newcommand{\E}{\mathbb{E}}
\renewcommand{\H}{\mathsf{H}}
\newcommand{\Q}{\mathsf{Q}}
\newcommand{\X}{\mathsf{X}}
\newcommand{\G}{\mathsf{G}}
\newcommand{\pair}[1]{\ensuremath{\left\langle #1 \right\rangle}}
\renewcommand{\Re}{\ensuremath{\mathbb{R}}}
\renewcommand{\prob}{\ensuremath{\mathbb{P}}}
\newcommand{\deriv}[2]{\ensuremath{\frac{\partial #1}{\partial #2}}}

\newcommand{\norm}[1]{\ensuremath{\left\| #1 \right\|}}

\newcommand{\braces}[1]{\ensuremath{\left\{ #1 \right\}}}
\newcommand{\parenth}[1]{\ensuremath{\left( #1 \right)}}

\newtheorem{definition}{Definition}
\newtheorem{theorem}{Theorem}
\newtheorem{remark}{Remark}

\begin{document}
\begin{frontmatter}

\title{Uncertainty propagation of stochastic hybrid systems: a case study for types of jump\thanksref{footnoteinfo}} 

\thanks[footnoteinfo]{This research has been supported in part by AFOSR under the grant FA9550-23-1-0400.}

\author[First]{Tejaswi K. C.} 
\author[Second]{William Clark} 
\author[Third]{Taeyoung Lee}

\address[First]{Mechanical and Aerospace Engineering, George Washington University, DC 20052 (email: kctejaswi999@gwu.edu)}
\address[Second]{Mathematics, Ohio University, Athens, OH 45701 (email: clarkw3@ohio.edu)}
\address[Third]{Mechanical and Aerospace Engineering, George Washington University, DC 20052 (email: tylee@gwu.edu)}

\begin{abstract}                
    Stochastic hybrid systems are dynamic systems that undergo both random continuous-time flows and random discrete jumps.
    Depending on how randomness is introduced into the continuous dynamics, discrete transitions, or both, stochastic hybrid systems exhibit distinct characteristics. 
    This paper investigates the role of uncertainties in the interplay between continuous flows and discrete jumps by studying probability density propagation. 
    Specifically, we formulate stochastic Koopman/Frobenius-Perron operators for three types of one-dimensional stochastic hybrid systems to uncover their unique dynamic characteristics and verify them using Monte Carlo simulations.
\end{abstract}


\end{frontmatter}

\section{Introduction}

Hybrid systems are dynamical systems that exhibit both continuous-time behaviors and discrete jumps.
Due to their capability to describe complex behaviors involving the interplay between continuous dynamics and discrete events in a unified manner, they have been actively investigated across a wide range of fields in dynamical system theory and control system design.

In particular, various formulations of stochastic hybrid systems have been proposed.
The General Stochastic Hybrid System (GSHS) is one of the most comprehensive frameworks for incorporating uncertainties in both continuous and discrete dynamics, as presented by \cite{bujorianu2006toward}.
The evolution of the continuous state is governed by a set of stochastic differential equations associated with each discrete mode.
The discrete jump process is triggered deterministically by entering a guard set or spontaneously according to a Poisson process.
The state after each jump is described by a stochastic kernel.
With its flexibility to incorporate uncertainties in every component of hybrid systems, GSHS has been widely applied in areas such as air traffic control~\citep{blom2009rare}, neuron modeling~\citep{pakdaman2010fluid}, and communication networks~\citep{hespanha2004stochastic}.

Uncertainty propagation examines how the distribution of the initial state evolves over time, and it is useful because the propagated density completely characterizes the statistical properties of the state at any given time.
To address this, the Frobenius-Perron operator was formulated by~\cite{oprea2023study} to propagate uncertainties for deterministic hybrid systems.
For GSHS with spontaneous jumps, spectral methods are applied to solve the associated hybrid Fokker-Planck equation in~\citep{wanleea19, wanleesjads22}.

However, depending on how randomness is incorporated, multiple types of stochastic hybrid systems arise, each with distinct dynamic characteristics.
For example, consider the two phases of a bouncing ball’s trajectory: the flight phase between bounces and the impact phase when the ball contacts the ground.
While the flight phase and the ball’s velocity after impact may be perturbed by various sources, there is no doubt that the ball will bounce upon hitting the ground.
This illustrates the case where random continuous-time flows represented by a stochastic differential equation are integrated with deterministic jumps and random resets, with other combinations also possible.

The objective of this paper is to investigate the interplay between continuous and discrete dynamics under different forms of randomness.
Specifically, we consider a one-dimensional hybrid system of three types: the deterministic case, random continuous flow with deterministic jumps, and random continuous flow with random spontaneous jumps.
For each case, we present the stochastic Koopman operator and the stochastic Frobenius-Perron operator, with the latter providing the Fokker-Planck equation to propagate uncertainties.
While we do not discuss it explicitly, the third case readily specializes to deterministic continuous flow combined with random jumps.
We demonstrate that the Koopman operator and the Frobenius-Perron operator yield distinct formulations for each case of stochastic hybrid systems, as verified by numerical studies involving Monte Carlo simulations.
In summary, the main contribution of this paper is the formulation of the Koopman and Frobenius-Perron operators for multiple types of stochastic hybrid systems, clarifying the role of uncertainties.

This paper is organized as follows.
Stochastic hybrid systems are introduced in \Cref{sec:SHS}, and the corresponding Koopman and Frobenius-Perron operators are formulated for each type of hybrid system in \Cref{sec:UP}, followed by numerical examples in \Cref{sec:NE}.

\section{Stochastic Hybrid Systems}\label{sec:SHS}

In this section, we formulate general stochastic hybrid systems (GSHS), which are uncertain dynamic systems that may undergo continuous-time evolution and discrete jumps. 
Then, we present the stochastic Koopman operator and the stochastic Frobenius-Perron operator. 

\subsection{General Stochastic Hybrid Systems}

The GSHS is defined by a tuple $\mathcal{H} = \left(\H,\mathcal{F}, \mu, X, h,\lambda,\G,K\right)$ described as follows:
\begin{itemize}
    \item Let $\Q=\{1,2,\ldots, n_q\}$ be the set of discrete modes, and let $\X_q\subset \Re^{n_q}$ with $n_q>0$ be the set of the continuous state in the $q$-th mode. 
        The hybrid state space is the disjoint union given by
        \begin{align}
            \H=\bigsqcup_{q\in\Q} \X_q =\bigcup_{q\in\Q}\{(x, q) | x\in\X_q\}.
        \end{align}
	The hybrid state is denoted by $(x,q)\in \H$, where $q$ is an auxiliary index that indicates which mode the state $x$ is from. 
	\item The $\sigma$-algebra on $\H$ is 
    \begin{align}
        \mathcal{F} = \{ \bigsqcup_{q\in\Q} A_q \,|\, A_q \in \mathcal{B}(\X_q)\},\label{eqn:F_hyb}
    \end{align}
    where $\mathcal{B}(\X_q)$ is the Borel $\sigma$-algebra on $\X_q$. 
    This corresponds to the set of all $A\subset \H$, satisfying $A\cap \X_q\in\mathcal{B}(\X_q)$ for any $q\in\Q$.
	\item The measure $\mu:\mathcal{F}\rightarrow\Re_{\geq 0}$ is
    \begin{align}
        \mu (A) = \sum_{q\in\Q} \mu_q(A \cap \X_q),\label{eqn:mu_hyb}
    \end{align}
    for $A\in\mathcal{F}$, 
    where $\mu_q:\mathcal{B}(\X_q)\rightarrow\Re_{\geq 0}$ is the Lebesgue measure on $(\X_q, \mathcal{B}(\X_q))$. 
	\item The continuous-time evolution of the state at each mode is governed by the following stochastic differential equation:
		\begin{equation} \label{eqn:ContSDE}
			dx = X(x, q)dt + h(x,q)dW,
		\end{equation}
	where $X(\cdot, q):\X_q\rightarrow\mathbb{R}^{n_q}$ is the drift vector, and $h(\cdot, q):\X_q\rightarrow\mathbb{R}^{n_q\times n_w}$ is the weighting matrix of diffusion for each $q\in \Q$. 
    And $W$ is a standard $n_w$-dimensional Wiener process.
	\item The spontaneous discrete transition is triggered by a Poisson rate function, $\lambda:\H\rightarrow\Re_{\geq 0}$, i.e., the probability that a transition occurs from $(x,q)$ over a time interval $\Delta t$ is $\lambda(x,q)\Delta t + o(\Delta t)$.
    \item Let the guard $\G=\sqcup_{q\in\Q} \G_q$ with $\G_q\subset \X_q$. The forced discrete transition is triggered when the hybrid state $(x,q)$ enters the guard set corresponding to the current mode $q$.
    \item During each discrete transition, the hybrid state is reset according to the stochastic kernel $K:\H\times\mathcal{F}\rightarrow\Re_{\geq 0}$. Specifically, $K((x,q)^-,H^+)$ corresponds to the probability measure for the set of posterior hybrid states $\H^+\in\mathcal{F}$ after jump, given the priori hybrid state $(x,q)^-\in\H$.
        When there is no spontaneous jump, the domain of the kernel reduces to $\G\times \mathcal{F}$. 
        In any case, $K(\cdot, \H)=1$. 
\end{itemize}

Under the given formulation of $(\H, \mathcal{F}, \mu)$, the Lebesgue integral of an integrable function $f$ on $\H$ is denoted interchangeably by
\begin{align*}
    \int f(x) d\mu(x) = \int f(x) \mu(dx).
\end{align*}
While the first expression is preferred, the second is used when the integration variable should be specified explicitly, especially when the measure itself may depend on another variable.

We make the following assumptions on GSHS. 
It is assumed that the number of jumps is finite for any finite interval of time, avoiding the Zeno behavior, and multiple jumps do not occur simultaneously. 
Further, we assume that the stochastic kernel for the discrete transition allows the duality in the following sense.
While the kernel is formulated as the probability measure of the posterior state after a jump given a priori state, it can be considered as an operation on a measurable function.
More specifically, for $f:\H\rightarrow\Re$, define $Kf:\H\rightarrow\Re$ as
\begin{align}
    (Kf)(x,q) = \sum_{q'\in\Q} \int_{\X_{q'}} f(x', q') K(x,q, dx', q').\label{eqn:K_operator}
\end{align}
This implies that when $(x,q)$ is a hybrid state right before a jump, then $(Kf)(x,q)$ corresponds to the mean of $f$ with respect to the posterior state after the jump.
Also, let the pairing on $f,g:\H\rightarrow\Re$ be
\begin{align}
    \pair{g(x,q), f(x,q)} = \sum_{q\in\Q}\int_{\X_q} g(x,q)f(x,q) d\mu(x).\label{eqn:pair}
\end{align}
We assume that there exists a dual of the stochastic kernel, $K^*:\H\times\mathcal{F}\rightarrow\Re_{\geq 0}$ in the sense that
\begin{align}
    \pair{g(x,q), (Kf)(x,q)} = \pair{(K^*g)(x,q), f(x,q)}.\label{eqn:K*}
\end{align}
In other words, $K^*$ is the dual of $K$ with respect to the pairing \eqref{eqn:pair}, when $K$ is perceived as a transformation on a measurable function as given by \eqref{eqn:K_operator}.


This formulation incorporates uncertainties in every component of hybrid systems: the randomness in evolution of the continuous state is given by the stochastic differential equation \eqref{eqn:ContSDE}; the stochastic jumps are modeled by a Poisson process where the jump rate function is dependent on the hybrid system or by a deterministic jump triggered by a guard; and the state right after the jump is described by a probability measure. 

\subsection{Stochastic Koopman/Frobenius-Perron Operator}




The stochastic Koopman operator and the Frobenious-Perron operator allow for a comprehensive framework for understanding and controlling stochastic systems by connecting probabilistic, data-driven, and dynamical perspectives.

Let $X_t$ be a time-homogeneous stochastic process indexed by $t\geq 0$.
When the process is deterministic, the Koopman operator $\mathcal{K}_t:L^\infty(\H)\rightarrow L^\infty(\H)$ is defined as $\mathcal{K}_t f(x) = f(X_t)$ when $X_0=x$ for $f\in L^\infty(\H)$. 
It is generalized to the stochastic setting by taking the expected value as follows~\citep{wanner2022robust}. 
\begin{definition}\label{def:SK}
    The stochastic Koopman operator $\mathcal{K}_t:L^\infty(\H)\rightarrow L^\infty(\H)$ for $t>0$ is defined as
\begin{align}
    \mathcal{K}_t f(x) & = \E[f(X_t)| X_0 = x] \label{eqn:K0},
\end{align}
for $f\in L^\infty(\H)$.
If the following limit exists,
\begin{align}
    \mathcal{A} f(x) 
    & = \lim_{t\downarrow 0} \frac{\mathcal{K}_{t}f(x) - f(x)}{t}, \label{eqn:AA}
\end{align}
then we state that $f$ belongs to the domain of $\mathcal{A}$, and $\mathcal{A}$ is the infinitesimal generator of $\mathcal{K}_t$. 
\end{definition}

In the literature of stochastic differential equations, \eqref{eqn:K0} is also referred to as transition operator semigroup, and several properties have been shown~\citep{oksendal2003stochastic}.
First, the Koopman operator is a semigroup, as given by
\begin{align}
    \mathcal{K}_{t+s}=\mathcal{K}_t\circ \mathcal{K}_s, \label{eqn:semigroup}
\end{align}
for any $t,s\geq 0$.
And the generator and the Koopman operator commute, i.e., $\mathcal{A}\mathcal{K}_t= \mathcal{K}_t\mathcal{A}$. 
Let $u(t,x) = \mathcal{K}_t f(x)$ for $f\in L^\infty(\H)$. 
If $f$ belongs to the domain of $\mathcal{A}$, then so does $u(t,x)$, and further, 
\begin{align}
    \deriv{u(t,x)}{t} & = \mathcal{A}u(t,x).\label{eqn:u_dot_Au}
\end{align}
By solving this differential equation, the evolution of the stochastic Koopman operator can be computed. 

In the deterministic setting, the Frobenius-Perron operator is considered as a pull-back of the differential form along the flow map~\citep{oprea2023study}.
However, such elegant formulation does not apply to the presented stochastic flow.
Instead, we formulate it as the dual of the stochastic Koopman operator as follows.
\begin{definition}
    Consider the stochastic Koopman operator as presented at \Cref{def:SK}.
    The corresponding stochastic Frobenius-Perron operator $\mathcal{P}_t: L^1(\H)\rightarrow L^1(\H)$ is defined such that
\begin{align}
    \pair{ g(x), \mathcal{K}_t f(x) } =  \pair{\mathcal{P}_t g(x), f(x)}, \label{eqn:adj_K_FP}
\end{align}
for $f\in L^{\infty}(\H)$ and $g\in L^1(\H)$. 
\end{definition}

From the duality given by \eqref{eqn:adj_K_FP}, we can show several properties of the stochastic Frobenius-Perron operator similar with \eqref{eqn:semigroup} and \eqref{eqn:u_dot_Au} of the Koopman operator.
First, using \eqref{eqn:semigroup} and \eqref{eqn:adj_K_FP}, we can show the semigroup property of the stochastic Frobenius-Perron operator, i.e.,
\begin{align}
    \mathcal{P}_{t+s} = \mathcal{P}_t \circ \mathcal{P}_s,
\end{align}
for any $t,s\geq 0$.
Also, the generator of the stochastic Frobenius-Perron operator is the adjoint of the generator of the stochastic Koopman operator, and they commute: $\mathcal{P}_t\mathcal{A}^* = \mathcal{A}^*\mathcal{P}_t$.
Further, similar with \eqref{eqn:u_dot_Au}, the stochastic Frobenius-Perron operator can be computed by
\begin{align}
    \deriv{v(t,x)}{t} = \mathcal{A}^* v(t,x), \label{eqn:v_dot_A*v}
\end{align}
where $v(t,x)=\mathcal{P}_t g(x)$ for $g\in L^1(\H)$. 

One of the most important roles of the stochastic Frobenius-Perron operator is that it propagates the probability density function. 
Suppose that there is a density function $p:\Re\times\H\rightarrow\Re_{\geq 0}$ such that
\begin{align*}
    \mathbb{P}[ X_t \in A ] = \int_{A}  p(t, x) d\mu(x),
\end{align*}
for $A\in\mathcal{F}$.
From the law of total expectation, for any $f:\H\rightarrow\Re$, we have
\begin{align*}
    \E[f(X_t)] & = \E[\E[f(X_t)| X_0=x]],
\end{align*}
where the outer expectation is taken with respect to the initial state $x$ distributed by $p(0,x)$,
and the inner expectation corresponds to the stochastic Koopman operator. 
Thus,
\begin{align*} 
    \E[f(X_t)] 
                  & = \int_\H p(0,x) \mathcal{K}_t f(x) d\mu(x)\\
                  & = \int_\H (\mathcal{P}_t p(0,x)) f(x)d\mu(x).
\end{align*}
where the second equality is from the duality.
As this holds for any $f\in L^\infty(\H)$,
this implies that $f(X_t)$ is distributed by $ \mathcal{P}_t p(0,x)$, which corresponds to the density at $t$, or $p(t,x)$.
In short, the stochastic Frobenius-Perron operator propagates the probability density as 
\begin{align}
    p(t,x)=\mathcal{P}_tp(0,x).\label{eqn:FP0}
\end{align}
As such, \eqref{eqn:v_dot_A*v} corresponds to the Fokker-Planck equation when $v(0,x)$ is chosen as the initial distribution.

\section{Uncertainty Propagation}\label{sec:UP}

In the preceding section, we presented a general formulation of hybrid dynamical systems, which can be affected by uncertainties in all aspects of how the state evolves between jumps (SDE), how the jump is initiated (guard or Poisson rate), and how the state is reset after the jump (reset kernel). 
In particular, one of the interesting characteristics of hybrid system is the interplay between the continuous-time flow and the discrete jump. 
Depending on whether the jump is triggered deterministically by intersecting a guard, or randomly according to a Poisson process, the stochastic flow exhibits distinct characteristics over the jump.

The objective of this section is to investigate the effects of the jump types on the propagation of uncertainties though a case study.
Specifically, we formulate the stochastic Frobenius-Perron operator for each of the following three types of hybrid systems:
\begin{itemize}
    \item Deterministic continuous flow with deterministic jumps
    \item Stochastic continuous flow with deterministic jumps
    \item Stochastic continuous flow with stochastic jumps 
\end{itemize}
Throughout this section, we focus on a one-dimensional example to highlight the unique characteristics of each case explicitly.

\subsection{Deterministic Case}

The first case is for the deterministic continuous flow integrated with deterministic jump. 
Specifically, the hybrid state space and the measure are chosen as $\Q=\{1\}$, $\X=(-\infty, b]$ for $b\in\Re$, and $d\mu=dx$. 

The continuous flow is governed by
\begin{align}
    \dot x = X(x). \label{eqn:dx1}
\end{align}
where $x\in\X$ and $X:\X\rightarrow\Re$ is a vector field. 
This corresponds to a special case of \eqref{eqn:ContSDE} when $h(x) = 0$.

Next, for deterministic jumps, the guard is chosen as the boundary of $\X$, or
\begin{align}
    \G = \partial \X = \{b\}. \label{eqn:G1}
\end{align}
There is no stochastic jump triggered by a Poisson process, i.e., $\lambda = 0$.
The reset kernel that specifies the posterior state after jump is given by a Dirac measure. 
Specifically, for $a\in\X$, 
\begin{align}
    K(x^-, H^+) = \begin{cases}
        1 & \text{if $a\in H^+$},\\
        0 & \text{otherwise},
    \end{cases}\label{eqn:K1}
\end{align}
for any $x^-\in\X$ and $H^+\in\mathcal{F}$. 
This implies that the state jumps from $b$ to $a$ always in a deterministic manner. 

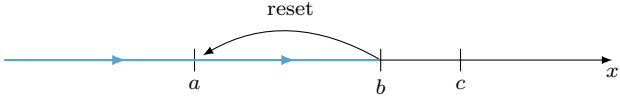
\begin{figure}
    \begin{tikzpicture}
        \footnotesize
        \draw[->] (0,0) -- (8,0) node[below] {$x$};
        \node (c) at (6,0) {};
        \node (a) at (2.5,0) {};
        \node[inner sep=0pt] (b) at ($(a)!0.7!(c)$) {};
        \draw (b |- 0, 0.15) -- ++(0,-0.3) node[below] {$b$};
        \draw (a |- 0,0.15) -- ++(0,-0.3) node[below] {$a$};
        \draw (c |- 0,0.15) -- ++(0,-0.3) node[below] {$c$};
        \draw[->] (b) to[out=150, in=30] node[midway, label=above: {reset}] {} (a);
        \draw[thick, myblue] (0,0) -- (b);
        \draw[->,myblue,thick] ($(a)!0.5!(b)$) -- ++(0.1,0);
        \draw[->,myblue,thick] (1.5,0) -- ++(0.1,0);
    \end{tikzpicture}
    \caption{1D hybrid system: the hybrid state space is $\X=(-\infty, b]$ (shaded by blue), and the guard is $\G=\partial\X = \{b\}$, at which the state jumps into $x=a<b$. The continuous vector field of $X(x)=-\gamma(x-c)$ for $\gamma>0$ and $c>b$ points to the right toward $x=c$.}\label{fig:case1}
\end{figure}

The resulting hybrid state space is illustrated at \Cref{fig:case1} for the special case when $X(x)=-\gamma(x-c)$ for $\gamma, c\in\Re$.

The Koopman operator and the Frobenius-Perron operator for this deterministic case are presented by \cite{oprea2023study}.
\begin{theorem}\label{thm:K1}
    Consider a hybrid system defined by \eqref{eqn:dx1}--\eqref{eqn:K1}. 
    The stochastic Koopman operator $u(t, x) = \mathcal{K}_t f(x)$ on $f:\X\rightarrow\Re$ satisfies
    \begin{align}
        \dfrac{\partial u(t,x)}{\partial t} & = \mathcal{A} u(t,x)\quad \mbox{for } x(t)\notin \{b\},\label{eqn:dotu1} \\
        u(t, a) & = u(t, b),\label{eqn:BCu1}
    \end{align}
    with the initial condition $u(0,x) = f(x)$, where the generator $\mathcal{A}:L^\infty(\X)\rightarrow L^\infty(X)$ is given by 
    \begin{align}
        \mathcal A u(t,x) &= X(x) \deriv{u(t,x)}{x} . \label{eqn:A1}
    \end{align}
\end{theorem}
\begin{pf}
    Let $\Phi_t:\X\rightarrow\X$ be the flow map for \eqref{eqn:dx1}. 
    Then, $u(t,x) = f(\Phi_t(x))$. From the chain rule, 
    \begin{align*}
        \deriv{u}{t} = \deriv{f(z)}{z}\bigg|_{z=\Phi_t(x)}\deriv{\Phi_t(x)}{t} = \deriv{u}{x} \dot x,
    \end{align*}
    which yields \eqref{eqn:dotu1} with \eqref{eqn:A1}.

    Next, suppose $x(t) = b$ for $t>0$, i.e., the trajectory intersects the guard at $t$.  
    According to the reset kernel, we have $\lim_{h\downarrow 0} \Phi_{t+h}(b) = a$.
    Thus, as $h\downarrow 0$, 
    \begin{align*}
        u(t+h, b) = \mathcal{K}_h (\mathcal{K}_t f(x)) = \mathcal{K}_h u(t, x)|_{x=b} = u(t,a),
    \end{align*}
    which provides \eqref{eqn:BCu1}.\hfill {} \qed
\end{pf}

Next, the Frobenius-Perron operator is constructed as follows. 
\begin{theorem}\label{thm:P1}
    Consider a hybrid system defined by \eqref{eqn:dx1}--\eqref{eqn:K1}. 
    The stochastic Frobenius-Perron operator $v(t, x) = \mathcal{P}_t g(x)$ on $g:\X\rightarrow\Re$ satisfies
    \begin{align}
        \dfrac{\partial v(t,x)}{\partial t} & = \mathcal{A}^* v(t,x)\quad \mbox{for } x(t)\notin \{a, b\},\label{eqn:dotv1} \\
        X(b) v(t, b)  & =   X(a)v(t, a^+) - X(a) v(t, a^-).\label{eqn:BCv1}
    \end{align}
    with the initial condition $v(0,x) = g(x)$, where the adjoint of the generator $\mathcal{A}^*:L^1(\X)\rightarrow L^1(X)$ is given by 
    \begin{align}
        \mathcal A^* v(t,x) &= -\deriv{(v(t,x)X(x))}{x}. \label{eqn:A*1}
    \end{align}
\end{theorem}
\begin{pf}
    Let $u(t,x)$ be the solution of \eqref{eqn:dotu1} and \eqref{eqn:BCu1} for the Koopman operator.
    From the duality given by \eqref{eqn:adj_K_FP}, we have
    \begin{align*}
        \pair{ v(t,x), \mathcal{K}_h u(t,x)} = \pair{ \mathcal{P}_h v(t,x), u(t,x)},
    \end{align*}
    for $h>0$. 
    Subtract $\pair{v(t,x), u(t,x)}$ from both side of the above equation to obtain
    \begin{align*}
        \pair{ v(t,x), \mathcal{K}_h u(t,x) - u(t,x)} = \pair{ \mathcal{P}_h v(t,x) - v(t,x), u(t,x)}.
    \end{align*}
    Divide both sides by $h$, and taking $h\downarrow 0$, 
    \begin{align}
        \pair{ v(t,x), \mathcal{A} u(t,x)} = \pair{ \mathcal{A}^* v(t,x), u(t,x)}.\label{eqn:uv_adj}
    \end{align}

    Substituting \eqref{eqn:A1}, the left hand side of the above is given by
    \begin{align*}
        \pair{v, \mathcal{A}u} = \int_{(-\infty, b]} v(t, x) X(x) \deriv{u(t,x)}{x} \, dx . 
    \end{align*}
    Due to the state reset at $x=a$, $v$ may not be continuously differentiable there. 
    As such, we split the domain of integration into $(-\infty, a-\epsilon]\cup(a+\epsilon, b]$ for $\epsilon\rightarrow 0$, and apply integration by parts to obtain
    \begin{align}
& \pair{v, \mathcal{A}u} =\int_{(-\infty, a^-]} v X \deriv{u}{x} dx  +\int_{(a^+, b]} vX \deriv{u}{x} dx\nonumber \\
& = vXu\bigg|_{a^-} + vXu\bigg|_{b} - vXu\bigg|_{a^+} - \int_{(-\infty, b]} \deriv{vX}{x} u dx. \label{eqn:thm2_0}
    \end{align}
    Since $u(t, a^+)=u(t, a^-)=u(t, b)$ from \eqref{eqn:BCu1}, the first three terms of the above reduce to
    \begin{align}
        u(t,b) \{ X(a)v(t,a^-) + X(b) v(t,b) - X(a)v(t, a^+) \}.\label{eqn:thm2_1}
    \end{align}
    We put together \eqref{eqn:uv_adj}, \eqref{eqn:thm2_0}, and \eqref{eqn:thm2_1} to obtain
 \eqref{eqn:BCv1} and \eqref{eqn:A*1}.\hfill {} \qed
\end{pf}

In the work by \cite{oprea2023study}, it is assumed that the reset state is another boundary of the hybrid state space. 
When applied to the presented one-dimensional example, this assumption implies that the hybrid state space is given by $\X=[a, b]$ instead of $\X=(-\infty, b]$. 
Despite this distinction in the problem formulation, we can show that the results of \cite{oprea2023study}, especially the concept of \textit{hybrid Jacobian}, yield a boundary condition consistent with \eqref{eqn:BCv1}. 

\begin{remark}
    In the development of \eqref{eqn:BCv1}, it is assumed that the trajectory \textit{intersects} the guard $x=b$. 
    This implies that the continuous flow is pointing toward the guard in the close neighborhood, i.e., $X(b) > 0$. 
    See \Cref{fig:case1}.
    When $X(b)\leq 0$, then there is no jump occurring at $x=b$, and the hybrid system reduces to continuous-time dynamics, where the boundary conditions \eqref{eqn:BCu1} and \eqref{eqn:BCv1} disappear. 
\end{remark}

\begin{remark}\label{rem:BC1}
    The boundary condition \eqref{eqn:BCv1} can be interpreted as the conservation of the total density. 
    Suppose $v(t,x)$ is the probability density of the state propagated by \eqref{eqn:dotv1}.
    The time rate of change of the total density is
    \begin{align*}
        \int_\X \mathcal{A}^* v(t,x) dx & = \int_{-\infty}^{a^-} - \deriv{vX}{x}dx + \int_{a^+}^b -\deriv{vX}{x} dx \\
                                                    & = - vX\big|_{a^-} + vX\big|_{a^+} - vX\big|_{b},
    \end{align*}
    which reduces to zero if \eqref{eqn:BCv1} is satisfied. 
    As such, \eqref{eqn:BCv1} can be interpreted as the balance between the rate of outflux at $b$ and the rate of influx at $a$. 
\end{remark}

\subsection{SDE with Deterministic Jump}

The next case is for the stochastic continuous flow integrated with deterministic jumps.
Compared with the first case, this incorporates randomness in the continuous-time evolution. 
The hybrid state space and the measure are same as the first case: $\Q=\{1\}$, $\X=(-\infty, b]$ for $b\in\Re$, and $d\mu=dx$. 
But, the continuous flow is governed by the following stochastic differential equation:
\begin{align}
    d x = X(x) dt + h(x) dW, \label{eqn:dx2}
\end{align}
where $x\in\Re$ and $X,h:\X\rightarrow\Re$, i.e., the special case of \eqref{eqn:ContSDE} for $\X\subset\Re$ and $\mathrm{dim}(\Q)=1$.

The deterministic jumps are identical to the first case, and they are defined by
\begin{align}
    \G & = \partial \X = \{b\}, \label{eqn:G2} \\
    \lambda & = 0,\\
    K(x^-, H^+) & = \begin{cases}
        1 & \text{if $a\in H^+$},\\
        0 & \text{otherwise},
    \end{cases}\label{eqn:K2}
\end{align}
for any $x^-\in\X$ and $H^+\in\mathcal{F}$. 

This second case incorporates an interesting interplay between a stochastic continuous flow and a deterministic jump. 
Similar with the first case, we develop the stochastic Koopman operator and the Frobenius-Perron operator as follows. 

\begin{theorem}\label{thm:K2}
    Consider a hybrid system defined by \eqref{eqn:dx2}--\eqref{eqn:K2}. 
    The stochastic Koopman operator $u(t, x) = \mathcal{K}_t f(x)$ on $f:\X\rightarrow\Re$ satisfies
    \begin{align}
        \dfrac{\partial u(t,x)}{\partial t} & = \mathcal{A} u(t,x)\quad \mbox{for } x(t)\notin \{b\},\label{eqn:dotu2} \\
        u(t, a) & = u(t, b),\label{eqn:BCu2}
    \end{align}
    with the initial condition $u(0,x) = f(x)$, where the generator $\mathcal{A}:L^\infty(\X)\rightarrow L^\infty(X)$ is given by 
    \begin{align}
        \mathcal A u(t,x) &= X(x) \deriv{u(t,x)}{x} + H(x)\frac{\partial^2 u(t,x)}{\partial x^2}. \label{eqn:A2}
    \end{align}
    with $H(x)=\frac{1}{2}h^2(x)$.
\end{theorem}
\begin{pf}
    The expression for the generator of \eqref{eqn:A2} can be derived by Ito's lemma, and it is well known in the literature, such as~\cite{oksendal2003stochastic}.
    The derivation of \eqref{eqn:BCu2} is identical to the proof of \Cref{thm:K1}. \hfill {} \qed
\end{pf}

\begin{theorem}\label{thm:P2}
    Consider a hybrid system defined by \eqref{eqn:dx2}--\eqref{eqn:K2}. 
    The stochastic Frobenius-Perron operator $v(t, x) = \mathcal{P}_t g(x)$ on $g:\X\rightarrow\Re$ satisfies
    \begin{gather}
        \dfrac{\partial v(t,x)}{\partial t}  = \mathcal{A}^* v(t,x)\quad \mbox{for } x(t)\notin \{a, b\},\label{eqn:dotv2} \\
        H(a) ( v(t, a^+) - v(t, a^-) ) = 0 ,\label{eqn:BCHa2} \\
        H(b) v(t, b) = 0,  \label{eqn:BCHb2} \\
        I (t, b) = I(t, a^+) - I(t, a^-),\label{eqn:BCI2}
    \end{gather}
    with the initial condition $v(0,x) = g(x)$, where the adjoint of the generator $\mathcal{A}^*:L^1(\X)\rightarrow L^1(X)$ is given by 
    \begin{align}
        \mathcal A^* v(t,x) &= -\deriv{v(t,x)X(x)}{x} + \frac{\partial^2 H(x)v(t,x)}{\partial x^2} = -\deriv{I(t,x)}{x}, \label{eqn:A*2}
    \end{align}
    and $I:\Re\times\X\rightarrow\Re$ is
    \begin{align}
        I(t,x) = v(t,x) X(x) - \deriv{H(x)v(t,x)}{x}. \label{eqn:I}
    \end{align}
\end{theorem}
\begin{pf}
    Let $u(t,x)$ be the solution of \eqref{eqn:dotu2} and \eqref{eqn:BCu2} for the Koopman operator.
    From the duality given by \eqref{eqn:adj_K_FP}, we have \eqref{eqn:uv_adj} satisfied as before.

    Substituting \eqref{eqn:A2}, the left hand side of \eqref{eqn:uv_adj} is given by
    \begin{align*}
        \pair{v, \mathcal{A}u} = \int_{(-\infty, b]} vX\deriv{u}{x} + vH\frac{\partial^2 u}{\partial x^2}\, dx.
    \end{align*}
    Similar with \eqref{eqn:thm2_0}, we split the domain of integration into $(-\infty, a-\epsilon]\cup(a+\epsilon, b]$ for $\epsilon\rightarrow 0$, and apply integration by parts to the last term to obtain
    \begin{align*}
\pair{v, \mathcal{A}u} 
& = vH\deriv{u}{x}\bigg|_{a^-} + vH\deriv{u}{x}\bigg|_{b} - vH\deriv{u}{x}\bigg|_{a^+} \\
& \quad + \int_{(-\infty, b]} I \deriv{u}{x} dx.
    \end{align*}
    The integration by parts is applied again for the last term, 
    \begin{align*}
        \pair{v, \mathcal{A}u} 
& = vH\deriv{u}{x}\bigg|_{a^-} + vH\deriv{u}{x}\bigg|_{b} - vH\deriv{u}{x}\bigg|_{a^+} \\
& \quad + I u\bigg|_{a^-} + Iu \bigg|_{b} - I u \bigg|_{a^+} 
 - \int_{(-\infty, b]} \deriv{I}{x} u dx,
    \end{align*}
    which is identical to $\pair{\mathcal{A}^* v, u}$ from \eqref{eqn:uv_adj}. 
    Since $u(t, a^+)=u(t, a^-)=u(t, b)$ from \eqref{eqn:BCu2} and $\deriv{u(t,a^+)}{x} = \deriv{u(t,a^-)}{x}$, we obtain the boundary conditions \eqref{eqn:BCHa2}--\eqref{eqn:BCI2} and the adjoint \eqref{eqn:A*2}.\hfill {} \qed
\end{pf}

\begin{remark}\label{rem:BC2}
    Similar with the discussion presented in \Cref{rem:BC1}, the boundary condition \eqref{eqn:BCI2} can be interpreted as the conservation of the total density. 
    The time rate of change of the total density is
    \begin{align*}
        \int_\X \mathcal{A}^* v(t,x) dx & = \int_{-\infty}^{a^-} - \deriv{I}{x}dx + \int_{a^+}^b -\deriv{I}{x} dx \\
                                                    & = - I\big|_{a^-} + I\big|_{a^+} - I\big|_{b},
    \end{align*}
    which reduces to zero if \eqref{eqn:BCI2} is satisfied. 
    Next, \eqref{eqn:BCHa2} corresponds to the continuity of the density at $a$, and \eqref{eqn:BCHb2} can be interpreted as the fact that the density is \textit{absorbed} at the guard $x=b$. 
\end{remark}

\begin{remark}
    While developed under district assumption, the presented case of SDE with deterministic jumps encloses the deterministic hybrid system as a special case. 
    By setting $H=0$, one can show that \Cref{thm:K2} and \Cref{thm:P2} reduce to \Cref{thm:K1} and \Cref{thm:P1}, respectively. 
\end{remark}

\subsection{SDE with Poisson Process}

The last case is for the stochastic continuous flow integrated with stochastic jumps.
Compared with the first two cases, this incorporates randomness in both of the continuous-time flow and discrete jump.
As there is no guard, the hybrid state space is expanded to $\X=\Re$. 
The set of modes and the measure are same as before: $\Q=\{1\}$,  and $d\mu=dx$. 
Same as \eqref{eqn:dx2}, the continuous flow is governed by the following stochastic differential equation:
\begin{align}
    d x = X(x) dt + h(x) dW, \label{eqn:dx3}
\end{align}
where $x\in\Re$ and $X,h:\X\rightarrow\Re$.

There is no guard triggering a jump in a deterministic way, i.e., $\G=\emptyset$. 
Instead there is a Poisson process with the rate function $\lambda:\X\rightarrow \Re_{\geq 0}$.
The kernel specifying the state after the jump is identical to the previous two cases, and it is given by
\begin{align}
    K(x^-, H^+) & = \begin{cases}
        1 & \text{if $a\in H^+$},\\
        0 & \text{otherwise},
    \end{cases}\label{eqn:K3}
\end{align}
for any $x^-\in\X$ and $H^+\in\mathcal{F}$. 
In this section, it is convenient to interpret the right hand side of \eqref{eqn:K3} as a Dirac measure, $\delta_a:\mathcal{F}\rightarrow\{0,1\}$, i.e.,
\begin{align}
    K(x^-, H^+) = \delta_{a} (H^+), \label{eqn:K3_dirac}
\end{align}
such that the action of the kernel on a function, namely \eqref{eqn:K_operator} is simplified into
\begin{align}
    (Kf) (x) = \int_\X f(x') d\delta_a( x' ) = f(a). \label{eqn:K_operator3}
\end{align}
Using the Dirac distribution (or Dirac delta function) $\delta_D:\Re\rightarrow\Re$, \eqref{eqn:K_operator3} can be further rearranged into
\begin{align}
    (Kf) (x) = \int_\X f(x') \delta_D(x'-a) dx'= f(a). \label{eqn:K_D_operator3}
\end{align}

\begin{theorem}\label{thm:K3}
    Consider a hybrid system defined by \eqref{eqn:dx3}--\eqref{eqn:K3}. 
    The stochastic Koopman operator $u(t, x) = \mathcal{K}_t f(x)$ on $f:\X\rightarrow\Re$ satisfies
    \begin{align}
        \dfrac{\partial u(t,x)}{\partial t} & = \mathcal{A} u(t,x),
    \end{align}
    with the initial condition $u(0,x) = f(x)$, where the generator $\mathcal{A}:L^\infty(\X)\rightarrow L^\infty(X)$ is given by $\mathcal{A} = \mathcal{A}_c+\mathcal{A}_d$ with
    \begin{align}
        \mathcal A_c u(t,x) &= X(x) \deriv{u(t,x)}{x} + H(x)\frac{\partial^2 u(t,x)}{\partial x^2}, \label{eqn:Ac3}\\
        \mathcal A_d u(t,x) &= \lambda(x)( (Ku)(t,x)-u(t,x)), \label{eqn:Ad3} \\
         &= \lambda(x)( u(t,a) -u(t,x)), \label{eqn:Ad3_dirac}
    \end{align}
    where $H(x)=\frac{1}{2}h^2(x)$.
\end{theorem}
\begin{pf}
    The expression for the generator of \eqref{eqn:Ac3} for the continuous flow is identical to \Cref{thm:K2}. 
    Next, let $U$ be the event that a jump occurs over the interval $[0,t]$ as $t\rightarrow 0$.
    From the formulation of the Poisson process, the probability for $U$ is $\prob[U] = t\lambda(x) +\mathcal{O}(t^2)$.
    After the jump, the posterior hybrid state $(x^+)$ is governed by the kernel $K$.
    Thus, 
    \begin{align}
        \E[u(t,X_t)|X_0=x, U] &= \int_{x^+\in\X} K(x, dx^+) u(t,x^+)\nonumber \\
          & = (Ku)(t,x),\label{eqn:EgU}
    \end{align}
    with $u(0,x)=f(x)$, 
    where the last equality is from \eqref{eqn:K_operator}.
    Next, the probability that no jump occurs is $\prob[U^c] = 1-t\lambda(x) + \mathcal{O}(t^2)$.
    Since the state is unchanged, the mean is given by
    \begin{align}
        \E[u(t,X_t)|X_0=x, U^c] = f(x) = u(0,x).\label{eqn:EgUc}
    \end{align}
    From the law of total expectation, we have
    \begin{align*}
    \E[u] = t\lambda \E[u|U] + (1-t\lambda)\E[u|U^c] +\mathcal{O}[t^2].
    \end{align*}
    Substituting \eqref{eqn:EgU}, \eqref{eqn:EgUc} and taking the derivatives with respect to $t$ at $t=0$, we obtain the generator of stochastic jumps as given by \eqref{eqn:Ad3}.

    Next, consider the specific kernel given by \eqref{eqn:K3}.
    From \eqref{eqn:K_operator3}, we have
    \begin{align*}
        (Ku) (t,x) & = u(t,a).
    \end{align*}
    Substituting this into \eqref{eqn:Ad3}, we obtain \eqref{eqn:Ad3_dirac}.  \hfill {} \qed
\end{pf}

\begin{theorem}\label{thm:P3}
    Consider a hybrid system defined by \eqref{eqn:dx3}--\eqref{eqn:K3}. 
    The stochastic Frobenius-Perron operator $v(t, x) = \mathcal{P}_t g(x)$ on $g:\X\rightarrow\Re$ satisfies
    \begin{gather}
        \dfrac{\partial v(t,x)}{\partial t}  = \mathcal{A}^* v(t,x),\label{eqn:dotv3}
    \end{gather}
    with the initial condition $v(0,x) = g(x)$, where the adjoint of the generator $\mathcal{A}^*:L^1(\X)\rightarrow L^1(X)$ is given by  $\mathcal{A}^*=\mathcal{A}_c^* + \mathcal{A}_d^*$ with
    \begin{align}
        \mathcal A_c^* v(t,x) &= -\deriv{v(t,x)X(x)}{x} + \frac{\partial^2 H(x)v(t,x)}{\partial x^2}, \label{eqn:Ac*3}\\
        \mathcal A_d^* v(t,x) &= (K^*(\lambda v))(t,x) - \lambda(x)v(t, x), \label{eqn:Ad*3} \\
                              & =\delta_D(x-a)\int_\X \lambda(y) v(t, y) dy - \lambda(x) v(t,x) \label{eqn:Ad*3_dirac}.
    \end{align}
    where $\delta_D:\Re\rightarrow\Re$ is the Dirac distribution. 
\end{theorem}
\begin{pf}
    The expression for the generator in \eqref{eqn:Ac*3} is identical to \Cref{thm:P2}. 
    Next, we construct $\mathcal{A}^*_d$ as follows. We have
    \begin{align*}
        \pair{ v, \mathcal{A}_d u} 
        & = \pair{v, \lambda (Ku) - \lambda u} \quad \mbox{from \eqref{eqn:Ad3}}\\
        & = \pair{\lambda v, (Ku) - u} \quad \mbox{as $\lambda$ is a scalar}\\
        & = \pair{K^*(\lambda v) - \lambda v, u},
    \end{align*}
    where the last equality is from \eqref{eqn:K*} and the fact that the pairing is linear. 
    This shows \eqref{eqn:Ad*3}. 

    Next, we show \eqref{eqn:Ad*3_dirac}.
    Using \eqref{eqn:K_D_operator3}, the left hand side of \eqref{eqn:K*} is
    \begin{align}
        \pair{ g(x), (Kf)(x)} & = f (a) \int_\X g(x) dx.\label{eqn:K*3_left}
    \end{align}
    We claim that the dual of \eqref{eqn:K3} is given by 
    \begin{align}
        K^*(x, dy) = \delta_D(x-a) dy. \label{eqn:K*3}
    \end{align}
    To verify \eqref{eqn:K*3}, we check the right hand side of \eqref{eqn:K*} as
    \begin{align*}
        \pair{ (K^*g)(x) ,f(x) } & = \int_\X \left( \int_\X g(y)\delta_D(x-a) dy\right)  f(x) dx\\
    & =  \int_\X g(y) dy \int_\X \delta_D(x-a) f(x) dx,
    \end{align*}
    which reduces to \eqref{eqn:K*3_left}. 
    Thus \eqref{eqn:K*3} is the dual of \eqref{eqn:K3}. 

    Substituting \eqref{eqn:K*3} into \eqref{eqn:Ad*3}, we obtain \eqref{eqn:Ad*3_dirac}.\hfill {} \qed
    \end{pf}

\begin{remark}
    The time rate of change of the total density is 
    \begin{align*}
        \int_\X \mathcal{A}^*v\, dx & =
        \int_\X \mathcal{A}_c^*v\, dx + \int_\X \mathcal{A}_d^*v\, dx,
    \end{align*}
    where the first term reduces to zero by following the discussion of \Cref{rem:BC2} under the assumption that $v$ is compactly supported in $\X$. 
    Next, the second term is
    \begin{align*}
        \int_\X \mathcal{A}_d^*v\, dx & = \pair{1, \mathcal{A}_d^*v} = \pair{\mathcal{A}_d 1, v}
                                      =\pair{ \lambda ((K1)-1), v},
    \end{align*}
    where we have used the duality and \eqref{eqn:Ad3}.
    But,  $K1(x) = \int_\X K(x, dy) dy = K(x, \X) = 1$ as $K$ is a probability kernel, which implies that the above reduces to zero. 
    As such, the total density is preserved under the flow of \eqref{eqn:dotv3}.
\end{remark}
\begin{remark}
    The generator \eqref{eqn:Ad3} and the adjoint \eqref{eqn:Ad*3} for the jump is not specific to the kernel $K$ defined at \eqref{eqn:K3}, and they can be applied to any other arbitrary kernel.
\end{remark}
\begin{remark}
    By setting $H=0$, \Cref{thm:K3} and \Cref{thm:P3} are applied to the case of hybrid systems where a deterministic continuous flow is integrated with spontaneous random jump, i.e., \eqref{eqn:dx1} with \eqref{eqn:K3}. 
\end{remark}

\section{Numerical Examples}\label{sec:NE}

In this final section, we simulate three different types of hybrid systems for the previously discussed one-dimensional case. 
We formulate finite-difference schemes for the governing equations of the Frobenius-Perron operator by discretizing both the spatial and temporal variables.
These schemes are solved using appropriate methods, taking the corresponding boundary conditions into account. 
To verify the characteristics of the obtained solutions, Monte Carlo simulations are constructed with a large number of particles, approximating the flow of the density function. 

\subsection{Deterministic Case}

\begin{figure}[h!]
	\centering
	\subfigure[Density evolution over time]{
		\includegraphics[width=0.475\linewidth]{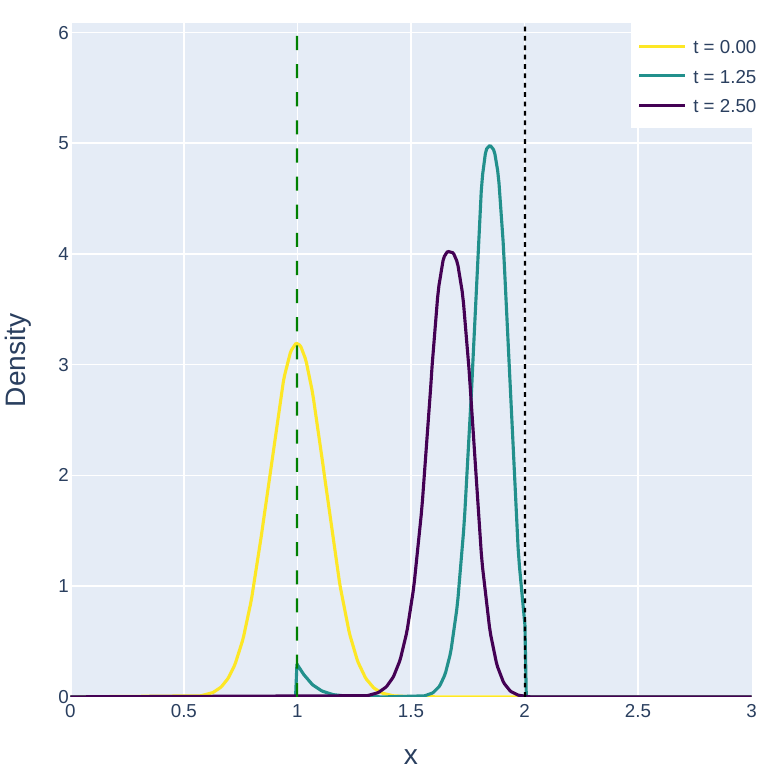}}
	\centering
	\subfigure[Monte carlo simulation]{
		\includegraphics[width=0.475\linewidth]{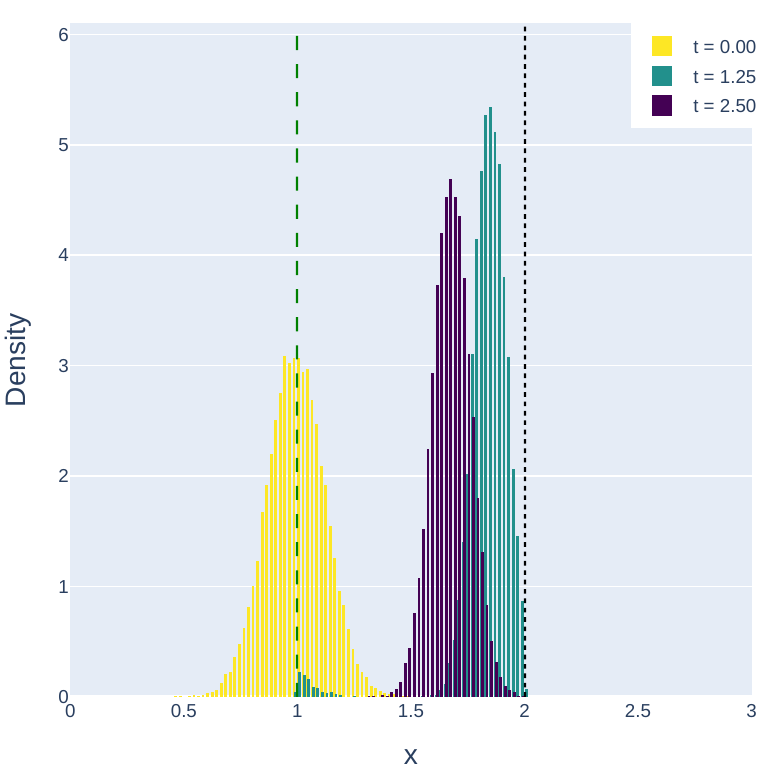}}
	\caption{Deterministic flow ($ H = 0 $)}
	\label{fig:det_case}
\end{figure}

Consider the hybrid system defined by \eqref{eqn:dx1}--\eqref{eqn:K1} with the specific vector field $ X(x) = -\gamma (x - c) $ as demonstrated in \Cref{fig:case1}. 
We use the parameters $ a = 1.0,\ b = 2.0,\ c = 3.0 $ and $ \gamma = 1.0 $. 

The density function $ v(t, x) $ evolves according to \eqref{eqn:dotv1} along with the jump condition at the boundary \eqref{eqn:BCv1}. 
This is a partial differential equation, specifically a nonlinear advection equation in the conservative form. 
Due to the expected discontinuity in the solution, it is challenging to implement finite difference schemes without introducing artificial dissipation. 
Therefore, we employ the MUSCL scheme~\citep{van1974towards}, a finite volume method with linear reconstruction and a \textit{minmod} flux limiter to avoid undesired oscillations near discontinuities.

For these illustrations, the initial condition is taken to be a normal distribution:
\begin{align*}
	v(0, x) = \mathcal{N}(a, \sigma^2)
\end{align*}
with $ \sigma = 0.125 $. 
The solution is computed at $ t = \braces{0, \cdots, t_k, t_{k+1}, \cdots} $ by advancing in time via an implicit scheme, which provides better accuracy compared to explicit methods.
Given $ v(t_{k}, x) $ in the spatial domain, the system of equations for the next time step is expressed as $ \varphi(v(t_{k+1}, x)) = 0 $ and solved using Newton's method. 

With these specific parameters values and initial condition, the Gaussian distribution travels rightward toward the boundary at $ x = b $.
Since the advection velocity $ |-\gamma(x-c)| $ is greater for particles further away from $ x = c $, the density is compressed leading to a higher peak and narrower spread near $ x = b $.
After reaching the boundary, particles are reset back to $ x = a $, causing this pattern to repeat periodically over time as observed in \Cref{fig:det_case}.(a).
Furthermore, the Monte Carlo simulation in \Cref{fig:det_case}.(b) confirms this solution since it started sampling from $ v(0, x) $ at $ t = 0 $.

\subsection{SDE with Deterministic Jump}

\begin{figure}[h!]
	\centering
	\subfigure[Density evolution over time]{
		\includegraphics[width=0.475\linewidth]{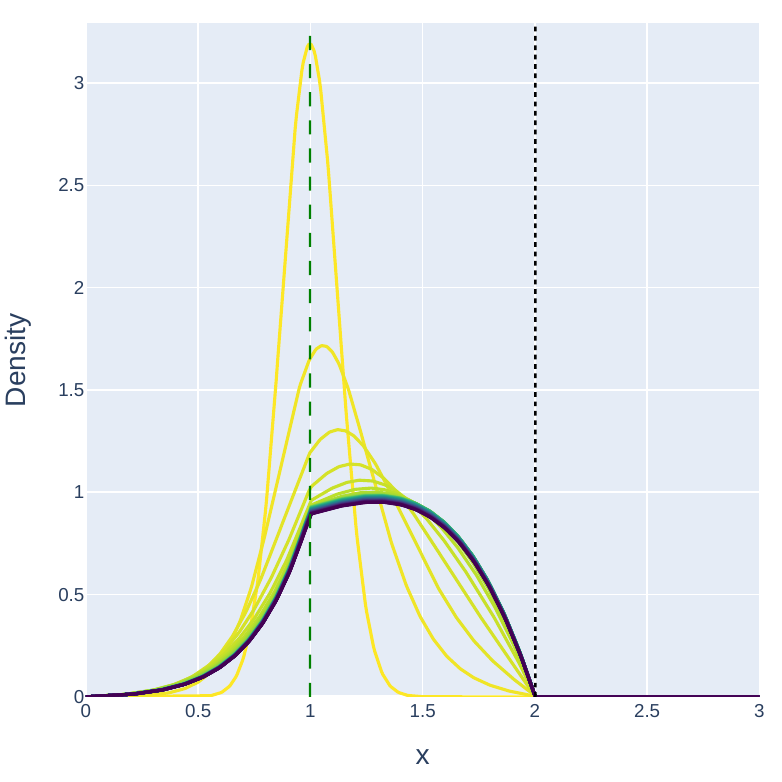}}
	\centering
	\subfigure[Monte carlo simulation]{
		\includegraphics[width=0.475\linewidth]{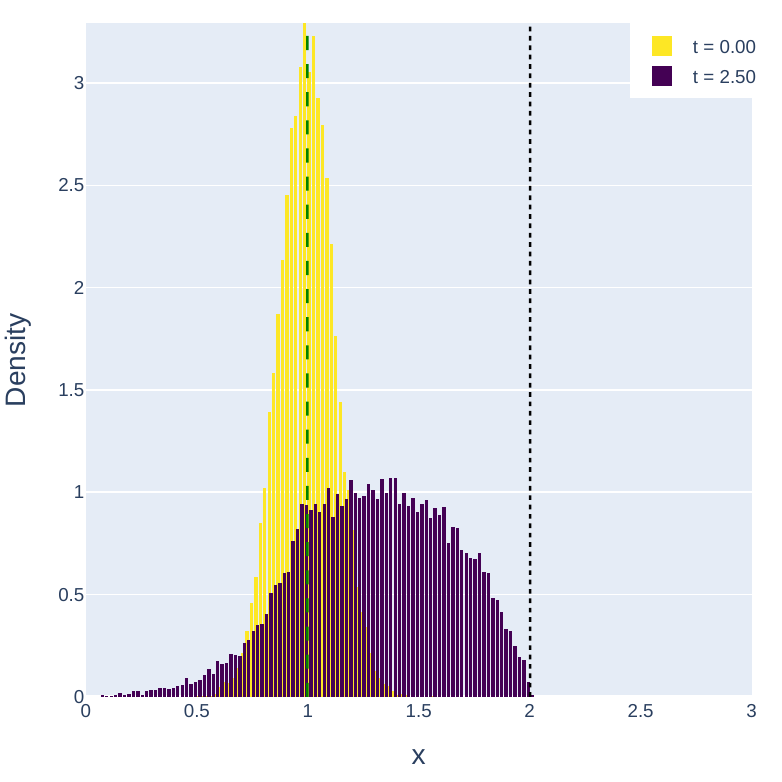}}
	\caption{Stochastic flow ($ H = 0.5 $) with deterministic jump}
	\label{fig:sde_det_jump_max}
\end{figure}

Next, we consider the hybrid system in \eqref{eqn:dx2}--\eqref{eqn:K2} which introduces additional noise into the continuous flow. 
The jump at the boundary $ x = b $ remains deterministic. 
Due to the stochastic component, the generator in \eqref{eqn:A*2} is augmented with a diffusion term that depends on $ H(x) $.
We apply Godunov's scheme~\citep{laney1998computational}, a classical upwind finite volume method for advection-diffusion equations, combined with an implicit formulation for time-stepping.

\begin{figure}[h!]
	\centering
	\subfigure[Density evolution over time]{
		\includegraphics[width=0.475\linewidth]{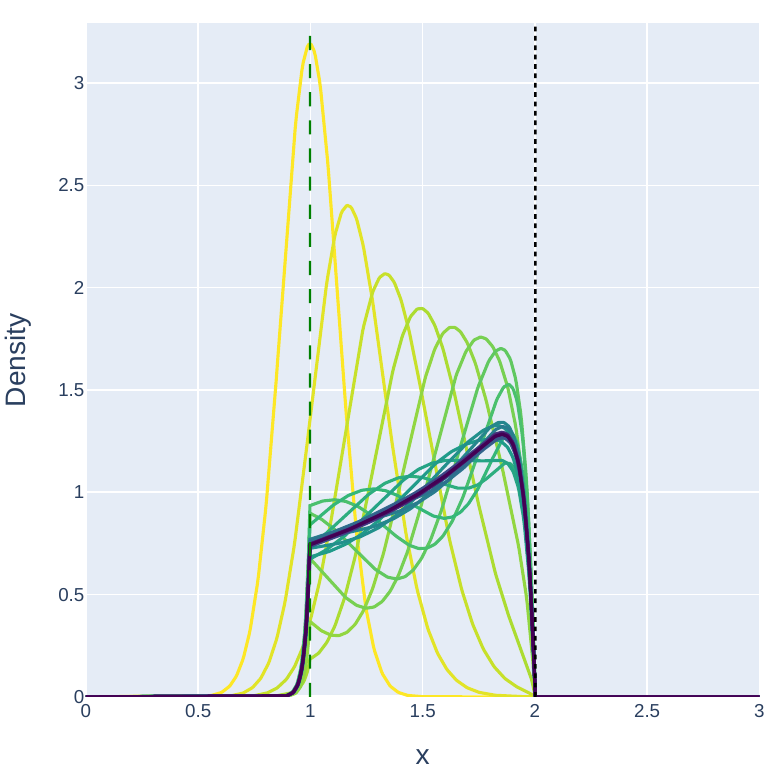}}
	\centering
	\subfigure[Monte carlo simulation]{
		\includegraphics[width=0.475\linewidth]{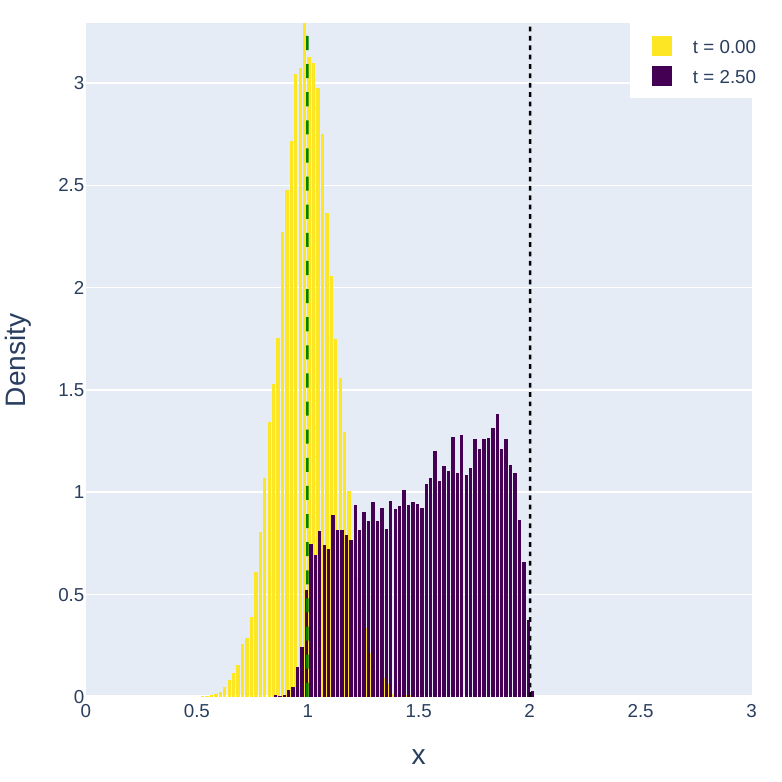}}
	\caption{Stochastic flow ($ H = 0.05 $) with deterministic jump}
	\label{fig:sde_det_jump}
\end{figure}

For non-zero stochasticity, the boundary conditions in \eqref{eqn:BCHa2}--\eqref{eqn:BCI2} simplify to $ v(t, a^-) = v(t, a^+),\ v(t, b) = 0 $ and 
\begin{align*}
	\deriv{H(x)v(t,x)}{x}\bigg|_{a^+} - \deriv{H(x)v(t,x)}{x}\bigg|_{a^-} = \deriv{H(x)v(t,x)}{x}\bigg|_{b}.
\end{align*}

We select parameters and an initial condition $ v(0, x) $ similar to those used in the previous case, to apply them to the governing equations in \eqref{eqn:dotv2}--\eqref{eqn:BCI2}. 
Two different diffusion cases are analyzed, $ H(x) \equiv 0.5 $ in \Cref{fig:sde_det_jump_max} and $ H(x) \equiv 0.05 $ in \Cref{fig:sde_det_jump}, over a simulation time horizon of $ T = 2.5 $ seconds.
It is interesting to observe the balance between the Gaussian distribution moving rightward towards $ x = c $ and the dissipation caused by the stochastic term. 
With a larger value of $ H $, even though particles reset to $ x = a $, there is a longer tail extending toward the negative $ x $-axis. 
Meanwhile, with a smaller $ H $, the density peak repeatedly resets at $ x = b $ and converges to the dark purple curve shown in \Cref{fig:sde_det_jump}.

In short, the diffusion in the continuous flow results in a stationary distribution in the limit, as opposed to the persistent jumps and resets observed in \Cref{fig:det_case}.
As discussed in \Cref{rem:BC2}, the density reduces to zero at $x=b$, and it is continuous at $x=a$. 
These boundary conditions are also distinct from \Cref{fig:det_case}, which are discontinuous and generally non-zero at those points. 

\subsection{SDE with Poisson Process}

Lastly, stochastic jumps are introduced to this continuous flow with noise. 
The jump itself is triggered by a Poisson process rather than a fixed guard set for the hybrid system in \eqref{eqn:dx3}--\eqref{eqn:K3}. 
Additional terms are included in \eqref{eqn:dotv3}--\eqref{eqn:Ad*3_dirac} from the discrete generator, $ \mathcal A_d^* v(t,x) $, specifically a decay term dependent on $ \lambda(x) $. 
Only within the spatial grid near $ x = a $, the first term in $ \mathcal A_d^* v(t,x) $ is non-zero.
Consequently, we have
\begin{align*}
\deriv{H(x)v(t,x)}{x}\bigg|_{a^+} - \deriv{H(x)v(t,x)}{x}\bigg|_{a^-} + \int_\X \lambda(y) v(t, y) dy = 0
\end{align*}
together with the continuity condition $ v(t, a^-) = v(t, a^+) $.

\begin{figure}[h!]
	\centering
	\subfigure[Density evolution over time]{
		\includegraphics[width=0.475\linewidth]{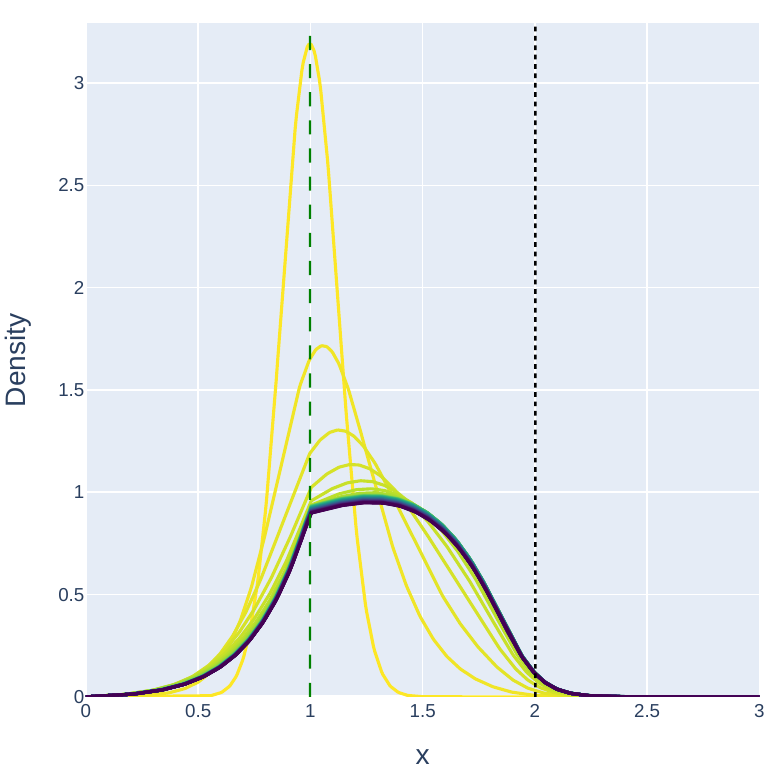}}
	\centering
	\subfigure[Monte carlo simulation]{
		\includegraphics[width=0.475\linewidth]{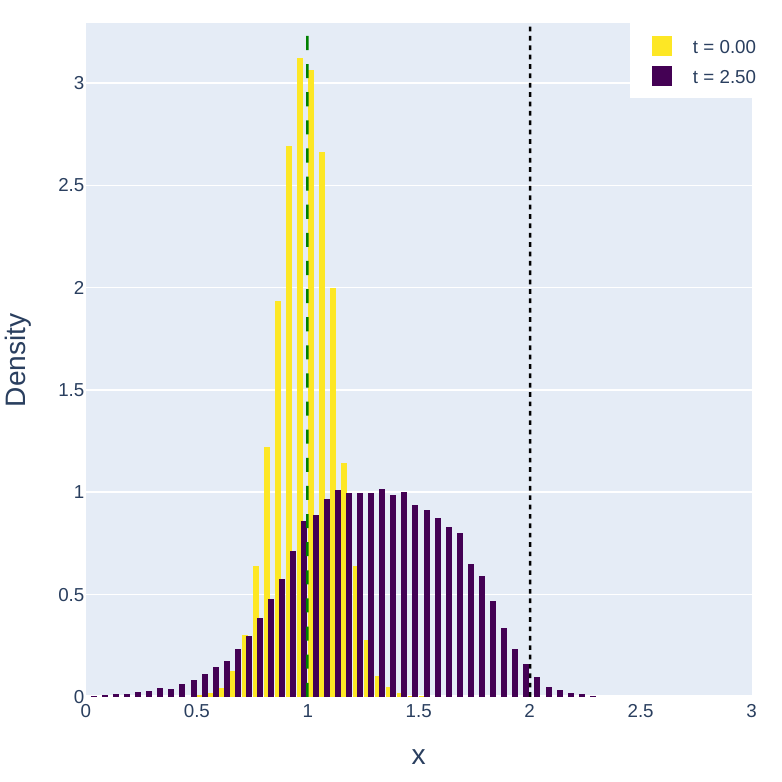}}
	\caption{Stochastic flow ($ H = 0.5 $) with Poisson process}
	\label{fig:sde_pois_jump_max}
\end{figure}

Using the same parameters and initial condition as before, two cases of $ H(x) $ are simulated.
We consider a continuous rate function given by
\begin{align*}
    \lambda(x) = 
    \begin{cases}
        0 & \text{if } x - b < - \epsilon,\\
        \frac{\lambda_{m}}{2} \parenth{1 + \sin\parenth{\frac{\pi}{2 x_{t}} (x - b)}} & \text{if } \norm{x - b} \le \epsilon,\\
        \lambda_{m} & \text{if } x - b > \epsilon,
    \end{cases}
\end{align*}
where a spatial threshold parameter $ \epsilon = 0.25 $ defines the boundary for reset, and the maximum rate is set to $ \lambda_m = 100.0 $.

\begin{figure}[h!]
	\centering
	\subfigure[Density evolution over time]{
		\includegraphics[width=0.475\linewidth]{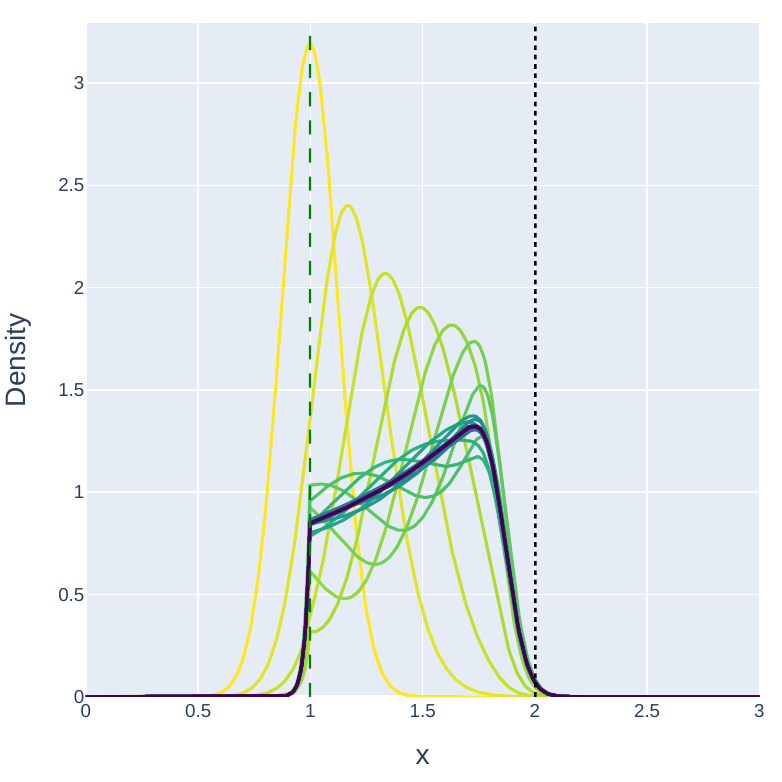}}
	\centering
	\subfigure[Monte carlo simulation]{
		\includegraphics[width=0.475\linewidth]{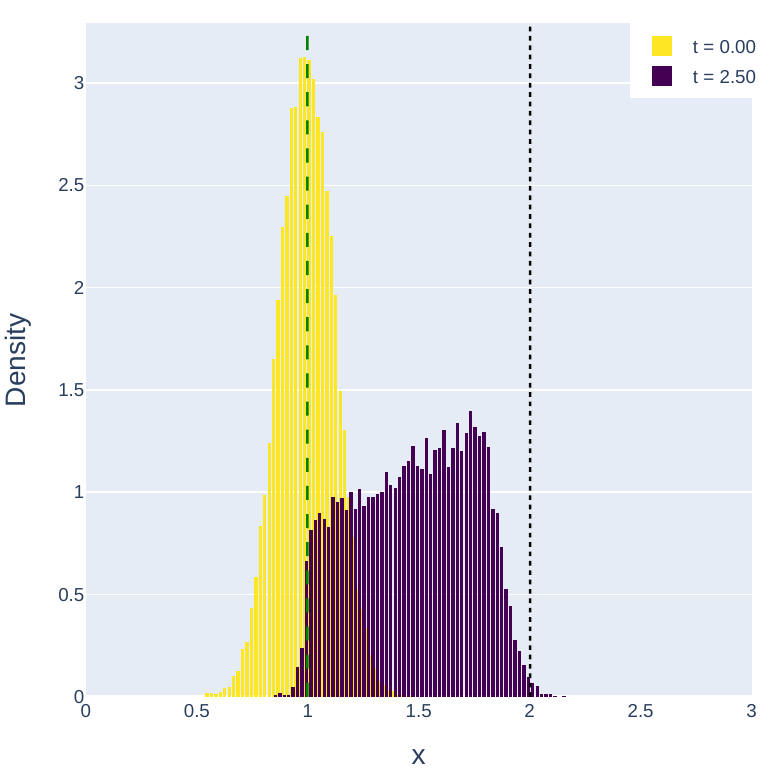}}
	\caption{Stochastic flow ($ H = 0.05 $) with Poisson process}
	\label{fig:sde_pois_jump}
\end{figure}

The corresponding results are illustrated at \Cref{fig:sde_pois_jump,fig:sde_pois_jump_max} for the two levels of diffusion.
Comparing with \Cref{fig:sde_det_jump,fig:sde_det_jump_max} for deterministic jumps, there are non-zero densities beyond the guard when $x>b$, which decays to zero as $x$ increases. 
This is no surprising as the jump is probabilistic: some trajectories may not jump at $x=b$. 
The response near $x=a$ is similar with the deterministic case.


\section{Conclusions}

This paper investigates the role of uncertainties in the interplay between the continuous flow and the discrete jump of stochastic hybrid systems. 
We formulate the stochastic Koopman operator and the stochastic Frobenius-Perron operator for three types of stochastic hybrid systems, where the latter is used to propagate uncertainties. 
It is demonstrated that depending on how randomness is incorporated into hybrid systems, the Frobenius-Perror operator exhibits distinct form and boundary conditions, and as such, the propagated density has unique characteristics. 
Future work includes generalizing these results to higher-dimensional manifolds.

\bibliography{ref} 

\begin{thebibliography}{11}
\providecommand{\natexlab}[1]{#1}
\providecommand{\url}[1]{\texttt{#1}}
\providecommand{\urlprefix}{URL }
\expandafter\ifx\csname urlstyle\endcsname\relax
  \providecommand{\doi}[1]{doi:\discretionary{}{}{}#1}\else
  \providecommand{\doi}{doi:\discretionary{}{}{}\begingroup
  \urlstyle{rm}\Url}\fi

\bibitem[{Blom et~al.(2009)Blom, Bakker, and Krystul}]{blom2009rare}
Blom, H.A., Bakker, G., and Krystul, J. (2009).
\newblock Rare event estimation for a large-scale stochastic hybrid system with
  air traffic application.
\newblock \emph{Rare event simulation using Monte Carlo methods}, 193--214.

\bibitem[{Bujorianu and Lygeros(2006)}]{bujorianu2006toward}
Bujorianu, M.L. and Lygeros, J. (2006).
\newblock Toward a general theory of stochastic hybrid systems.
\newblock In \emph{Stochastic hybrid systems}, 3--30. Springer.

\bibitem[{Hespanha(2004)}]{hespanha2004stochastic}
Hespanha, J.P. (2004).
\newblock Stochastic hybrid systems: Application to communication networks.
\newblock In \emph{International Workshop on Hybrid Systems: Computation and
  Control}, 387--401. Springer.

\bibitem[{Laney(1998)}]{laney1998computational}
Laney, C.B. (1998).
\newblock \emph{Computational gasdynamics}.
\newblock Cambridge university press.

\bibitem[{{\O}ksendal(2003)}]{oksendal2003stochastic}
{\O}ksendal, B. (2003).
\newblock \emph{Stochastic differential equations}.
\newblock Springer.

\bibitem[{Oprea et~al.(2023)Oprea, Shaw, Huq, Iwasaki, Kassabova, and
  Clark}]{oprea2023study}
Oprea, M., Shaw, A., Huq, R., Iwasaki, K., Kassabova, D., and Clark, W. (2023).
\newblock A study of the long-term behavior of hybrid systems with symmetries
  via reduction and the frobenius-perron operator.
\newblock \emph{arXiv preprint arXiv:2309.12569}.

\bibitem[{Pakdaman et~al.(2010)Pakdaman, Thieullen, and
  Wainrib}]{pakdaman2010fluid}
Pakdaman, K., Thieullen, M., and Wainrib, G. (2010).
\newblock Fluid limit theorems for stochastic hybrid systems with application
  to neuron models.
\newblock \emph{Advances in Applied Probability}, 42(3), 761--794.

\bibitem[{Van~Leer(1974)}]{van1974towards}
Van~Leer, B. (1974).
\newblock Towards the ultimate conservative difference scheme. {II}.
  {M}onotonicity and conservation combined in a second-order scheme.
\newblock \emph{Journal of computational physics}, 14(4), 361--370.

\bibitem[{Wang and Lee(2020)}]{wanleea19}
Wang, W. and Lee, T. (2020).
\newblock Spectral bayesian estimation for general stochastic hybrid systems.
\newblock \emph{Automatica}, 117.
\newblock \doi{10.1016/j.automatica.2020.108989}.

\bibitem[{Wang and Lee(2022)}]{wanleesjads22}
Wang, W. and Lee, T. (2022).
\newblock Uncertainty propagation for general stochastic hybrid systems on
  compact {Lie} groups.
\newblock \emph{SIAM Journal on Applied Dynamical Systems}, 21(3), 2215--2240.
\newblock \doi{10.1137/21M144147X}.

\bibitem[{Wanner and Mezic(2022)}]{wanner2022robust}
Wanner, M. and Mezic, I. (2022).
\newblock Robust approximation of the stochastic koopman operator.
\newblock \emph{SIAM Journal on Applied Dynamical Systems}, 21(3), 1930--1951.

\end{thebibliography}


\end{document}